\RequirePackage{fix-cm}
\documentclass[smallextended]{svjour3}       
\smartqed  
\usepackage{graphicx}
\usepackage{color}
\usepackage{subfigure}
\usepackage{url}
\usepackage{amsmath,amsfonts,amssymb}
\usepackage{array,multirow}
\usepackage{graphics}
\usepackage{float}
\usepackage{caption}
\usepackage{verbatim}
\usepackage{color}
\usepackage{rotating}
\usepackage{tabularx}
\usepackage{lscape}
\usepackage[table]{xcolor}
\usepackage{tikz}
\usepackage{pgfplots}
\usepackage{pgfplotstable}
\usepackage[normalem]{ulem}
\usetikzlibrary{arrows,%
                petri,%
                topaths,
 shapes,trees}%

\usepackage{algorithmic}

\usepackage{booktabs}

\def\tr#1#2#3{\langle#1,#2,#3\rangle}

\begin{document}

\title{Drone-assisted deliveries: new formulations for the Flying Sidekick Traveling Salesman Problem
}
%


\titlerunning{New formulations for the FSTSP}        

\author{Mauro Dell'Amico         \and
        Roberto Montemanni \and
        Stefano Novellani
}


\institute{Mauro Dell'Amico \at
              Dipartimento di Scienze e Metodi dell'Ingegneria (DISMI), Universit\`a di Modena e Reggio Emilia (UNIMORE)\\
            via Amendola 2, 42122 Reggio Emilia, Italy
           \and
           Roberto Montemanni \at
              Dipartimento di Scienze e Metodi dell'Ingegneria (DISMI), Universit\`a di Modena e Reggio Emilia (UNIMORE)\\
            via Amendola 2, 42122 Reggio Emilia, Italy
           \and Stefano Novellani \at
              Dipartimento di Scienze e Metodi dell'Ingegneria (DISMI), Universit\`a di Modena e Reggio Emilia (UNIMORE)\\
            via Amendola 2, 42122 Reggio Emilia, Italy \\
              Tel.: +39-0522-523537\\
              \email{stefano.novellani@unimore.it}           
}

\date{Received: 29/03/2019 / Accepted: date}

\maketitle

\begin{abstract}
In this paper we consider a problem related to deliveries assisted by an unmanned aerial vehicle, so-called drone. In particular we consider the Flying Sidekick Traveling Salesman Problem, where a truck and a drone cooperate to delivery parcels to customers minimizing the completion time. In the following we improve the formulation found in the related literature. We propose three-indexed and two-indexed formulations and a set of inequalities that can be implemented in a branch-and-cut fashion.
We could find the optimal solutions for most of the literature instances. Moreover, we consider two versions of the problem: one in which the drone is allowed to wait at the customers, as in the literature, and one where waiting is allowed only in flying mode. The solving methodologies are adapted to both versions. A comparison between the two versions is provided.
\end{abstract}



\section{Introduction}\label{sec:intro}

Once restricted to the military domain, {\em unmanned aerial vehicles} (UAV), also known as (aerial) drones, have received a widespread adoption in civil applications by
{the humanitarian sector and are currently object of great interest of the commercial sector}.
This is due to their capability of accomplishing otherwise impossible tasks, performing  activities more effectively, or to lead to cost saving solutions.
A comprehensive definition of drones that is the one provided by
Scott and Scott \cite{scott2017drone}: \textit{``Drones are devices which are capable of sustained flight, which do not have human on board, and are under sufficient control to perform useful functions"}.

During the last years, drones have interested researchers whose effort was put especially in improving the required technology:
\begin{itemize}
\item[(a)] the hardware, by reducing weigh, increasing battery duration, improving the charging system, improving safety, etc.,
\item[(b)]  the software, by designing better autonomous operations and guidance systems using  improved  GPS accuracy, localization techniques, obstacle detection and avoidance techniques, enhanced sensors, image processing, etc.,
\item[(c)] the safety and security elements, by adopting technologies to protect the flight against spoofing and hijacking of UAVs, etc.
\end{itemize}
A smaller effort has been done for the solution of optimization and operational problems related to the use of drones, although the number of paper in this direction has strongly increased in the very last years.

{
The paper is organized as follows: we start with Section \ref{sec:app} where we summarize the main
real-life applications that use drones. In Section \ref{sec:lit} we discuss upon the most relevant related literature, while in Section \ref{sec:prob} we formally describe the problem. Mathematical formulations and their implementations are described in Sections \ref{sec:ImprovementOnF1} and \ref{sec:2form}.
Extensive computational results comparing the models performances are presented in Section \ref{sec:results}. The last Section \ref{sec:conclusions} concludes the paper.
}

\section{Applications} \label{sec:app}
Early drones application arose in the military sector, where most of the optimization works are focused on the {\em path planning}, normally to avoid obstacles or radars, minimizing travel length or altitude, or to improve the usage of batteries. In this paper we are not interested in listing the military applications, but we address the interested reader to the following works: Bortoff \cite{bortoff2000path}, Richards et al. \cite{richards2002coordination}, Zheng et al. \cite{zheng2005evolutionary}, and Roberge et al. \cite{roberge2013comparison}.

Besides military applications, drones can be used in both humanitarian and commercial sectors.
The three sectors make use of drones in such a way that can be divided in two main applications: to collect and deliver information and to collect and deliver goods. On one hand we have the surveillance, monitoring, or covering, and on the other hand deliveries. In the following we describe a large set of applications adhering to this dichotomy.

\subsection{Gathering information: Surveillance, monitoring and covering}

The information gathering applications include surveillance, monitoring, and covering activities, where displacement of goods is not necessary, in which drones fly autonomously and monitor the environment with different sensors or cameras and communicate and exchange data and information with other drones or with a central station.

Drones can be used to optimize the coverage of an area (see, e.g., Shang et al. \cite{shang2014algorithm}) or following
specific targets that can also be moving ones (see, e.g., Zorbas et al. \cite{zorbas2013energy} and Di Puglia Pugliese et al. \cite{pugliese2016modelling}).
This is done by defining the position of drones, their path (see, e.g., Kashuba et al. \cite{kashuba2015optimization}), their number, etc.{, while
maximizing the coverage or}  minimizing the time between an appearance of an event and its covering, the total length, the service costs in energy, etc.
{In particular,  typical applications regard the following domains: traffic management, environmental monitoring, catastrophic events, remote locations surveillance, precision agriculture, building inspection, security surveillance, etc.}

\subsection{Moving goods: Deliveries}

This second set of applications is characterized by the fact that real goods need to be moved from one place to another. The main application is to provide faster, more cost efficient goods or parcel delivery, especially in the last-mile, but also mail delivery and quickly required medical items delivery have been object of interest. Indeed, in emergency and events management drones can be used when infrastructure are damaged or unreliable to access isolated regions to deliver needed goods in addition to surveillance. This point can match the medical items delivery: delivery of blood, medications, vaccines, defibrillators, insulin, oxygen, or other needed health-care items in location with difficult access due to poor infrastructure, remote areas, traffic congestion, inaccessible roads due to weather or disasters, or simply urgently needed.
Drone can transport water, food, medical supplies during a crisis or and event.

{Commercial applications arise when one needs to provide faster, more cost efficient goods or parcel delivery, especially in the last-mile. }
Several companies are investigating the use of drones for parcel delivery for e-commerce, such as Amazon \cite{amazon}, Alibaba \cite{alibaba}, and Alphabet \cite{alphabet}.
Amazon CEO announced Amazon's Prime Air, that uses a fleet of UAVs to deliver parcels from warehouses to customers \cite{amazon}.
Australian textbook distributor Zookal started testing drone parcel deliveries in Australia \cite{zookal}.
DHL Parcel operates an autonomous drone delivery system, to deliver medications and other urgent goods to one of the Germany's North See islands of Juist \cite{dhlisland}. In 2016 they also tested a delivery system in the Bavaria alpine region typically carrying either sporting goods or urgently needed medicines \cite{dhlalps}.
Chinese JD.com deploys drones to extend their delivery and logistics network. They consider to use drones in areas with complex terrain and poor infrastructure for last-mile delivery. They started in 2016 with four rural locations of China in the outskirts of Beijing and in the provinces of Jiangsu, Shaanxi, and Sichuan. They announced an agreement with the Shaanxi provincial government to build China's largest low altitude drone network, serving a 300 km area with drone stations and routes to delivery e-commerce parcels \cite{JD,JDshaanxi}. They will build 150 drone launch facilities in China by 2020 \cite{jd2020}.
Also Alphabet (Google) announced to enter the field with Project Wing \cite{alphabet}.
UPS and DPDgroup are also testing parcel delivery with drones \cite{ups,dpd}.
UPS and Zipline are working on a drone network to deliver vaccines and blood to 20 clinics in remote locations in Rwanda \cite{upsrwanda}. Only 20\% of Africans live within 2 km of a road that functions year-round, thus malaria medications, rabies vaccines, etc. can be carried and left with a parachute (see Scott and Scott \cite{scott2017drone}).
UPS and Workhorse are also testing drone and truck combo delivery \cite{upswork}.
Flirtey completed the first Federal Aviation Administration approved drone delivery in July 2015, when it delivered medical supplies to a health clinic in Wise, Virginia. The company has also started a partnership with 7-Eleven for home delivery \cite{flirtey}.
In the Netherlands and Sweden a prototype ambulance drone has also been tested for delivering defibrillators \cite{nederland,sweden}.
Matternet provides an on-demand delivery platform, an end-to-end solution integrating the Matternet's drones and stations. They provide their platform as a service to healthcare, e-commerce and logistics organizations. In particular, they transport medical items between health-care facilities in Switzerland. They perform drone deliver in Zurich. In this case Matternet and Mercedes-Benz have joined forces to create delivery solutions integrating vans and drones for siroop online shop \cite{mercedes}.
The United Arab Emirates say they plan to use unmanned aerial drones to deliver official documents and packages to its citizens as part of efforts to upgrade government services \cite{emirates}.

This is just an incomplete list of many of the cases that arose in the field in the last years. The large number make us suspect an increasing interest in the field in the future.

\section{Related literature} \label{sec:lit}
The FSTSP is a generalization of the (TSP) and the {\em vehicle routing problem} (VRP). A large body of literature has been dedicated to these problems; however, we are interested only in those problems that are highly related to the FSTSP, especially where drones and trucks are coupled and synchronized.\\

A problem that is conceptually related to the FSTSP  is the {\em Close enough traveling salesman problem} (see, e.g., Shuttleworth et al. \cite{shuttleworth2008advances}), that aims at finding the cheapest route for the truck without visiting every customer on his route, but only getting within a certain radius of each customer.
A similar problem and even more conceptually related is the {\em covering salesman problem} (see, e.g., Current and Schilling \cite{current1989covering}),
that aims at finding the cheapest tour such that all nodes that are not part of the route lie within a specified radius from a node on the truck route. This reminds the maximum radius that the drone can travel.

In the FSTSP, the vehicle that arrives first at the meeting point has to wait for the other, thus we can state that it lies in the class of problems that require synchronization between vehicles, in particular, among the categories defined by Drexl \cite{drexl2012synchronization}, the FSTSP can be considered under {\em movement synchronization en route}, where vehicles may join and separate multiple times along a route. The author claims that this class received little attention in literature. One of these problems is also one of the most related to our problem: the {\em Truck and Trailer Routing Problem} (TRRP) (see e.g., Chao \cite{chao2002tabu}), where two different types of vehicle can serve customers: trucks and trailers. Due to practical constraints (e.g. street size), some customers can only be served by a truck, while other customers can be served by a truck or a truck pulling a trailer. A truck is autonomous, while the trailers always need to be pulled by a truck. Parking places are used to decouple trucks from trailer when convenient to serve some subset of customers that can only be served by single truck. This resembles the idea of the drone truck decoupling.
A similar problem to the FSTSP has been studied by Lin \cite{lin2011vehicle}, where two types of delivery resources are used, vans and foot couriers, which allows for coordination. A heavy resource (a van) may carry both delivery items and one or more units of the lighter resource (foot couriers). Foot couriers can pick up and deliver items independently or travel with a van on its outbound and/or return leg, they can serve more than one customers, in contrary of what a drone can do. Moreover, foot couriers do not need to return to the same van.\\

Let us now consider works that solve problems related to the use of drones.
As said, the first applications of UAVs have been performed in the military sector, and the same can be said for the use of optimization techniques when optimizing operations for drones. In Sisson's thesis \cite{sisson1997applying} it is studied a tabu search coupled with Monte Carlo simulation to determine the minimum number of UAVs to cover a pre-selected target set based on stochastic survival probabilities that also incorporate the wind effects.
In Ryan's thesis \cite{ryan1998embedding} a tabu search within a discrete event simulation is applied to solve a multi TSP with time windows for UAVs, where it is required to attain a level of target coverage using a minimum number of vehicles. Weather and threat are considered.
A direct extension of this work is proposed in Ryan et al. \cite{ryan1998reactive}.
O'Rourke et al. \cite{o1999dynamic} consider dynamic routing of UAVs in operational use with the US Air Force. Dynamic components are wind and emerging targets. They model the problem as a VRP with time windows and use a tabu search to solve it.

Boone et al. \cite{boone2015enhanced} solve a multi TSP for drones, they firstly divide customers in clusters by using K-mean clustering, then they solve a TSP for each cluster with a nearest neighbor improved with 2-opt.
Dorling et al. \cite{dorling2017vehicle} study a multi trip VRP for drone delivery where the effect of batteries and weight on energy consumption is considered. They solve one version of the problem in which costs are minimized under a delivery time limit and a second version where delivery time is minimized subject to a budget constraint. The authors present MILP formulations and a simulated annealing algorithm.
Tseng et al. \cite{tseng2017flight} solve a modified TSP where an autonomous drone has to serve delivery points and may use charging stations for charging the battery. They consider wind uncertainty and variable speeds. The used algorithm is based on the Christofides one to solve the shortest spanning tree and then the minimum matching problem to obtain an Eulerian tour.

Location problems linked to drones are studied in the two following works. Scott and Scott \cite{scott2017drone} consider drone delivery models for healthcare. They consider delivery with trucks that leave a central depot to drone nests, and then from drone nests to delivery points. They use two objective functions: one minimizes the total time and the other minimizes the maximum weighted time for truck/drone delivery. Both depot and nests need to be located. No routing is considered but only direct trips.
Shavarani et al. \cite{shavarani2017application} studied the facility location problem for the optimization of drone delivery system evaluating the Amazon prime air case study in San Francisco. They locate both launching stations and recharge stations. No routing is considered but direct trips. The authors solve the problem with a genetic algorithm.

Campbell et al. \cite{campbell2018drone} present some drone arc routing problems, that are continuous optimization problems. They discretized them by approximating the curves with polygonal chains. The drone rural postman problem, in which drones have no capacity, is solved with iterative algorithms for the rural postman problem. They discuss also other problems where drones are capacitated or more than one.  \\

The following papers consider the coupled interventions of drones and trucks. The first paper we consider is slightly different from the others, indeed Savuran and Karakaya \cite{savuran2015route} study a problem where a truck follows a linear path, in the meanwhile a drone is launched from the truck and must return to the truck after performing a, let's say, open TSP to visit all the targets. They solve the problem with a genetic algorithm.
Ferrandez et al. \cite{ferrandez2016optimization} propose a work that couples drones and a truck. They have a set of delivery customers that are clusterized by using K-means clusterization. They thus solve the TSP among the centroid of each cluster that are the points where the truck stops to launch one or more drones. The TSP part is solved with a genetic algorithm. Boysen et al. \cite{boysen2018drone} consider a fixed truck route where the truck represents a loading platform for the drones. The truck and the drone can wait each other. They use at most 2 drones and three restriction of the problem: one where the drone must return to the launching point, one in which the drone can return up to the next vertex, and the last one where the drones can return in one of the following nodes. \\

In the following we consider problems where drones and vehicles work together to complete operations, in which both can accomplish tasks. In this first part we evaluate what Otto et al. \cite{ottooptimization} classify as {\em Drones and vehicles performing independent tasks}. Murray and Chu \cite{murray2015flying} propose the {\em parallel drone scheduling TSP} (PDSTSP), where one truck and a fleet of drones can serve customers only departing from the depot. In this case only customers within a certain range from the depot can be served via drones. The others are served by the truck. Solving this problem can provide good results when many customers lie close to the depot. They propose MILP formulations and simple greedy heuristics for both problems. Saleu et al. \cite{mbiadou2018iterative} propose a two step iterative heuristic based on dynamic programming for the same problem.
Ulmer and Thomas \cite{ulmer2018same} study a dynamic variant of the PDSTSP called the {\em Same-day delivery with heterogeneous fleets of drones and vehicles}, where requests arrive dynamically and they need to be allocated to drones or truck maximizing the number of  served customers. They solve the problem with an approximate dynamic programming known as parametric policy function approximation.

We consider, now, routing problems where trucks are equipped with drones and both vehicles can be used to deliver packages to customers. Otto et al. \cite{ottooptimization} classify these problems under the name {\em Drones and vehicles as synchronized working units}. In the case of the drone, the flight is called {\em sortie}: the drone start from its vehicle in a vertex of the network (launch), performs a delivery to a customer, and returns to the truck (rendezvous) in a vertex of the network. Some problems consider the launching and rendezvous times negligible, some other account for these; however, truck and drone must be synchronized and thus wait for each other. The objective function is normally to minimize the completion time.
Murray and Chu \cite{murray2015flying} define and study the {\em Flying sidekick traveling salesman problem} (FSTSP). In the FSTSP truck and drone can cooperate to serve customers: one drone can leave the truck at a vertex
and return to the truck at another vertex
after completing a delivery. Customers can be visited only once, but some customers can be visited only by the truck because their request cannot be fulfilled by the drone (for various reasons, such as capacity limitations, requirement of a signature, drone cannot safely land, etc.). Drones cannot return at the launching point.
Agatz et al. \cite{agatz2018optimization} solve the {\em TSP with Drone} (TSP-D), where one truck cooperates with one drone to make deliveries. The aim is to find the fastest method to serve customers with a truck or a drone that can leave the truck, serve a customer, and return to the truck in one of the vertices.
Each customer has to be visited at least once by one of the vehicles, but they can be visited more than once by the truck if it is convenient for drone launching and return. Launching and rendezvous of the same sortie can coincide.
Some vertices cannot be visited by drones. Endurance is unlimited and launching and rendezvous times are considered negligible.
The authors preset an ILP and propose route first-cluster second heuristics based on local search and dynamic programming. Bouman et al. \cite{bouman2018dynamic} solve the TSP-D with with dynamic programming. The drone is not slower than the truck. Instances with up to 20 customers are solved.
Ha et al. \cite{ha2015heuristic} also use the name TSP-D, albeit in this case nodes cannot be visited multiple times and launch and rendezvous of a sortie are not allowed to happen at the same vertex. The authors declare that the considered problem shared the FSTSP characteristics, but they called it TSP-D. To highlight the difference between the two problems we call this problem TSP-D*. The authors propose two heuristic algorithms: a route first-cluster second one and a cluster first-route second one.
Ponza's thesis \cite{ponza2016optimization} tackles the FSTSP proposing a modified MILP formulation with respect to the Murray and Chu's \cite{murray2015flying} one and solve it with a simulated annealing algorithm. The new formulation, with respect to the FSTPS one, among the other few differences, does not allow the drone to wait at customers nodes.

Wang et al. \cite{wang2017vehicle} define the {\em vehicle routing problem with drones} (VRPD), where a homogeneous fleet of trucks equipped with a not necessary unitary number of drones delivers parcels to customers. Drones can be launched from trucks at depot or at any customer vertex. Each drone must return to the same truck also if this happens at the same node where it has been launched. One drone can carry only one parcel. The authors imagine that the drone must travel along the street network.
They derive a number of worst case results based on the number of drones per truck and the difference between the drones speed and the trucks speed.
Poikonen et al. \cite{poikonen2017vehicle} extend the worst-case results considering different metrics for trucks and drones, considering limited drone batteries, and evaluating different objective functions. They also evaluate connections between the VRPD and other VRP variants.
Daknama and Kraus \cite{daknama2017vehicle} present and solve the {\em Vehicle routing with drones} where multiple vehicles and drones can be used for deliveries. They minimize the average delivery times instead of the completion time. No mathematical model is presented, but they solve the problem by first solving a multiple TSP heuristically and then introducing drones. Local search procedures are thus applied to improve the solution. \\

In Table \ref{tab:problems} we classify the truck and drone problems. In $Paper$ column we report the authors and the reference to the paper, in the second column we report the name of the problem as by the papers (but we include the renaming of TSP-D* to avoid confusion:
we recall that Agatz et al.  \cite{agatz2018optimization} and Ha et al. \cite{ha2015heuristic} use the same name for different problems
).
In $\#t$ and $\#d$ we report the considered number of trucks and drones, respectively. If trucks and drones are multiple we identify it with the letter $m$. In following columns we consider some features of the problems: the fact that launch and return points may coincide ($L=R$); the fact that trucks may visit vertices multiple times ($m-v$); the fact that some vertices can be visited only by trucks ($t-v$). In last column we provide a rough description of the method used to solve the proposed problems, if any. {Table \ref{tab:problems} is inspired by the tables proposed by Ponza \cite{ponza2016optimization}
that we adjusted referring only to the problems where drones and vehicles are considered as synchronized working units and by adding the fact that the studied problems allow the following points: (i) the possibility of having multiple visits; (ii) the possibility of having nodes to be visited only by the truck.}

{A recent survey by Otto et al. \cite{ottooptimization} provides a wide overview on  civil applications of drones, it gives an insight into optimization approaches used to solve operational problems, in particular where both drones and drones combined with other vehicles are considered, that is the class of problems addressed in this work.}

\begin{table}[htbp]
\tiny
  \centering

    \begin{tabular}{m{6em}lm{5em}rrrrrm{11em}}
    \toprule
    Paper & Name  & Formul. &
    {\#t} & {\#d} & {L=R} & {m-v} & {t-v} & Solving method \\
    \midrule
   Murray and Chu \cite{murray2015flying} & FSTSP & MILP 2-index & 1 & 1 &   &   & x & MILP, heuristic (nearest neighbour, savings, sweep + savings for UAV routes definition) \\

    Agatz et al. \cite{agatz2018optimization} & TSP-D & ILP set covering & 1 & 1 & x & x & x  & Heuristic (route first cluster second) \\

    Bouman et al. \cite{bouman2018dynamic} & TSP-D & none  & 1 & 1 & x & 
    x & x  & Dynamic programming \\

    Ha et al. \cite{ha2015heuristic} & TSP-D* & {MILP for clustering, set packing} & 1 & 1 &   &   &   & Heuristic (route first cluster second, cluster first route second) \\

    Ponza \cite{ponza2016optimization}& FSTSP & MILP, 2-index & 1 & 1 &  &   & x & MILP, saimulated annelaing \\

    Wang et al. \cite{wang2017vehicle} & VRPD & none & m & m & x &   &   & none \\

    Poikonen et al. \cite{poikonen2017vehicle} & VRPD & none & m & m & x &   &   & none \\

    Daknama and Kraus \cite{daknama2017vehicle} & VRD & none & m & m &  &   &   & Metaheuristic, Local search \\

    \bottomrule
    \end{tabular}%
  \caption{Classification of problems where drones  and  vehicles  are  synchronized working  units.}\label{tab:problems}%
\end{table}%

\section{Problem Description and Basic Mathematical Model} \label{sec:prob}

In this work we study the FSTSP defined by Murray and Chu \cite{murray2015flying}, that is to serve a set of customers $C = \{1,\dots,c\}$ with either a truck or a drone. The truck starts from depot 0 and returns to the final depot $c+1$, and is equipped with a flying drone that can be used in parallel to serve one customer at a time. The drone can perform a {\em sortie}, defined by a launching node (where the drone leaves the truck), a served customer, and a rendezvous node (where the drone returns to the truck), that must be different from the launching one. All customers of $C$ can be served by the truck, but only a subset $C' \subseteq C$ can be served by the drone with a sortie. The problem is built on digraph $G=(N,A)$, where the set $N = \{0,1,\dots,c+1\}$ represents all the nodes, while we define $N_0 = \{0,1,\dots,c\}$ and $N_+=\{1,\dots,c+1\}$. The set $A$ is the set of all the arcs $(i,j), i \in N_0, j \in N_+$, $i\neq j$. Each arc $(i,j)$ is associated with two non-negatives traveling times: $\tau_{ij}^T$ and $\tau_{ij}^D$, that represent the time for traveling that arc by the truck and by the drone, respectively. The travel time matrix of the drone and the truck may be different ($\tau_{ij}^T\lessgtr\tau_{ij}^D, (i,j)\in A$). Nodes 0 and $c+1$ represent the same physical point, the depot, and the traveling time between them is set to 0 to account for the case in which there is only one customer served by the drone directly from the depot.

{We assume that the capacity of the truck is large enough to serve all customers and that the drone performs only one delivery at the time, i.e., it leaves the truck, serves a customer and returns to the truck before possibly serving a new customer.}
A {sortie} is formally defined by a triplet $\langle i,j,k \rangle$, ($i\ne j\ne k$) where $i \in N_0$ is the launching node, $j \in C'$ the customer to serve, and $k \in N_+ $ the rendezvous node.

Serving times at customers for both drone and truck are negligible.
Service times for preparing the drone at launch and rendezvous are given by $\sigma^L$ and $\sigma^R$.
Drones have a battery limit (endurance) of $E$ time units, that limits drone use. {Rendezvous time $\sigma^R$ contributes to the endurance computation while $\sigma^L$ does not, since the drone lies on the truck when it is prepared for launch.}

Let $F$ be the set of all sorties that can be performed within the endurance $E$,  that means
$\tau_{ij}^D + \tau_{jk}^D + \sigma_R \le E$, $i \in N_0, j \in C', k \in N_+$.

Drones can be launched from trucks only when the truck is stopped at customers points or at the depot{, however, they cannot leave the depot before the truck starts its route.} In the meanwhile a truck can keep serving customers and the drone can return only at another point, that is not the launching one (if not the depot). This requires a certain synchronization: the vehicle (drone or truck) that arrives first at the meeting point has to wait for the other.
The objective is to minimize the completion time, that is the moment when the last vehicle arrives at the depot.

\subsection{Murray-Chu  Formulation} \label{sec:murray}

We report in this section the mathematical formulation proposed by Murray and Chu \cite{murray2015flying} (MC in the following) that is the basis for our work.
A binary variable $x_{ij}$ is created for each arc $(i,j) \in A: i \ne j$ and is set to 1 if the truck uses the corresponding arc, 0 otherwise. For each  sortie of the drone is used a binary variable $y_{ijk}, \langle i,j,k \rangle \in F$ that takes value one if the sortie
is performed, 0 otherwise.

The non-negative variables $t_i^T, i \in N$ and $t_i^D, i \in N$ represent the time of availability at node $i$
for the truck and the drone, respectively. For the starting depot we fix $t_0^T=t_0^D=0$.
  Note that the drone and the truck may leave uncoupled at the starting depot. The waiting of the truck (resp. of the drone) at a node $i$ is modeled as a delayed
availability and included in $t_i^T$ (resp. $t_i^D$).

The MC formulation uses two sets of auxiliary variables. A
variable $u_i, i \in N$, $1 \le u_i \le c+2$ is used to model  subtour elimination constraints for the track route, as in  Miller-Tucker-Zemlin {(see, e.g., Desrochers and Laporte \cite{desrochers1991improvements}).} Finally, binary variable $p_{ij}$ is set to 1 if customer $i\in C$ is visited before customer $j \in C$ $(j\ne i$), thus defining a total ordering among all the pair of customers. We set $p_{0j} = 1$, $j \in C$ to impose the starting of the route from the depot.

The resulting model for the MC formulation is described in the following, presenting one group of logically related constraints at a time.

The objective function \eqref{eq:om1_fobb}  minimizes the  arrival time at the depot of the truck, but due to next  constraints \eqref{eq:om1_droneTruck3} and \eqref{eq:om1_droneTruck4} this is equivalent to minimize $\max\{t_{c+1}^T,t_{c+1}^D\}$.
\allowdisplaybreaks
\begin{align}
\min  t_{c+1}^T\label{eq:om1_fobb}
\end{align}

\noindent
\emph{Customer covering} \\
The first constraints impose that all customers must be served either by the truck or by the drone:
\begin{align}
&\sum_{(i,j)\in A} x_{ij} +  \sum_{  \tr{i}{j}{k} \in F} y_{ijk} = 1 \quad  j \in C\label{eq:om1_ass}
\end{align}
\noindent\emph{Truck routing constraints}\\
Constraints \eqref{eq:om1_start} and \eqref{eq:om1_end} force the truck to depart from depot $0$ and to return to depot $c+1$ at the end of the trip, while
constraints \eqref{eq:om1_flow} guarantee the flow conservation at customers.\\
\allowdisplaybreaks
\begin{align}
&\sum_{\substack{j \in N_+}} x_{0j} = 1 \label{eq:om1_start}\\
&\sum_{\substack{i \in N_0}} x_{i,c+1} = 1 \label{eq:om1_end}\\
& \sum_{(i,j)\in A} x_{ij} = \sum_{(j,k)\in A} x_{jk} \quad j \in C \label{eq:om1_flow}
\end{align}
\noindent\emph{Single sortie leave/return}\\
Constraints \eqref{eq:om1_dep}--\eqref{eq:om1_arr} impose that at most one launch of drone and one
rendezvous with the truck is done  in any node.
The drone cannot start its sorties in $c+1$ or return to $0$.
\allowdisplaybreaks
\begin{align}
&\sum_{\tr{i}{j}{k}\in F} y_{ijk}  \le 1 \quad i \in N_0 \label{eq:om1_dep}\\
& \sum_{\tr{i}{j}{k}\in F}  y_{ijk}  \le 1 \quad k \in N_+ \label{eq:om1_arr}
\end{align}
\noindent\emph{$x$-$y$ coupling constraints}\\
In constraints \eqref{eq:om1_link} we impose that, if the triplet $\tr{i}{j}{k}$ is selected, then the truck must enter in node $i$ and in node $k$ to launch and collect the drone. Constraints \eqref{eq:om1_link0} is the equivalent of \eqref{eq:om1_link} when  the drone is launched from the depot.
\begin{align}
& 2y_{ijk} \le \sum_{(h,i)\in A} x_{hi} + \sum_{(l,k)\in A} x_{lk} \quad  \tr{i}{j}{k}\in F, i\in N_+ \label{eq:om1_link}\\
& y_{0jk} \le \sum_{(h,k)\in A} x_{hk} \quad \tr{0}{j}{k}\in F \label{eq:om1_link0}
\end{align}
\noindent\emph{Truck-drone timing constraints}\\
Constraints from \eqref{eq:om1_droneTruck1} to \eqref{eq:om1_droneTruck4}  ensure time synchronization between the truck and the drone at launching and return node.
Constraints \eqref{eq:om1_tt0} and \eqref{eq:om1_td0} set the starting time at depot of both drone and truck to zero.
\allowdisplaybreaks
\begin{align}
&t_i^D \ge t_i^T - M (1-\sum_{\tr{i}{j}{k}\in F} y_{ijk}) \quad i \in C\label{eq:om1_droneTruck1}\\
&t_i^D \le t_i^T + M (1-\sum_{\tr{i}{j}{k}\in F} y_{ijk}) \quad i \in C\label{eq:om1_droneTruck2}\\
&t_k^D \ge t_k^T - M (1-\sum_{\tr{i}{j}{k}\in F} y_{ijk}) \quad k \in N_+\label{eq:om1_droneTruck3}\\
&t_k^D \le t_k^T + M (1-\sum_{\tr{i}{j}{k}\in F} y_{ijk}) \quad k \in N_+\label{eq:om1_droneTruck4}\\
&t_0^T = 0\label{eq:om1_tt0}\\
&t_0^D = 0\label{eq:om1_td0}
\end{align}
\noindent\emph{Truck timing constraints}\\
Constraints \eqref{eq:om1_time5} state that if the arc $(h,k)\in A$ is used by the truck, then
the timing variables must be consistent with the travelling times and the launch and rendezvous times, if drone is used.
We will discuss these constraints in detail in  Section \ref{subsec:F1}.
\begin{align}
\begin{split}t_k^T \ge t_h^T + \tau_{hk}^T +  \sigma_R \sum_{\tr{i}{j}{k}\in F}  y_{ijk} + \qquad\qquad\qquad\qquad\qquad\qquad\\
\sigma_L\sum_{\substack{\tr{k}{l}{m}\in F\\ l\neq h}} y_{klm} - M(1- x_{hk})  \quad  (h,k)\in A
\end{split} \label{eq:om1_time5}
\end{align}
\noindent\emph{Drone timing constraints}\\
Constraints \eqref{eq:om1_droneTime6}-\eqref{eq:om1_droneTime7} impose consistency on the drone timing variables when a sortie $\tr{i}{j}{k}$ is selected.
\begin{align}
&t_j^D \ge t_i^D + \tau_{ij}^D - M (1-\sum_{\tr{i}{j}{k}\in F}y_{ijk} ) \quad (i,j)\in A, j \in C' \label{eq:om1_droneTime6}\\
&t_k^D \ge t_j^D + \tau_{jk}^D + \sigma_R - M (1-\sum_{\tr{i}{j}{k}\in F}y_{ijk} ) \quad (j,k)\in A, j \in C'\label{eq:om1_droneTime7}
\end{align}
Note that these constraints  allow $t^D_j$ to be greater than the arrival time at customer $j$, and consequently to insert a waiting time at $j$.\\

\noindent\emph{Drone battery endurance constraints}\\
Constraints \eqref{eq:om1_energy} are the battery endurance constraints, for each sortie.
\begin{align}
&t_k^D -(t_j^D -\tau_{ij}^{D}) \le E + M (1-y_{ijk} ) \quad \tr{i}{j}{k}\in F \label{eq:om1_energy}
\end{align}
We observe that inequalities \eqref{eq:om1_energy} do not consider the waiting at $j$ in the computation of the energy consumption, i.e., we can therefore suppose that the drone is allowed to wait on the ground.\\
\noindent\emph{Miller-Tucker-Zemlin subtour elimination}\\
Constraints \eqref{eq:om1_u_sub} are the classical Miller-Tucker-Zemlin constraints for the truck route, while  \eqref{eq:om1_u_sub_2} are the extension to sorties, and impose that node $i$ precedes node $k$ if sortie $\tr{i}{j}{k}$ is selected.\\
\allowdisplaybreaks
\begin{align}
&u_i - u_j + 1 \le (c+2)(1-x_{ij}) \quad (i,j)\in A, i\in N_+\label{eq:om1_u_sub}\\
&u_k - u_i \ge  1 - (c+2)(1-\sum_{\tr{i}{j}{k}\in F}y_{ijk}) \quad i \in C, k \in N_+, k\ne i\label{eq:om1_u_sub_2}
\end{align}

\noindent\emph{Node total ordering}\\
Constraints \eqref{eq:om1_p1}--\eqref{eq:om1_p0} induce variables $p_{ij}$ to define a total ordering of the nodes.
\begin{align}
&p_{ij} + p_{ji} = 1 \quad (i, j)\in A, i <j\label{eq:om1_p1}\\
&p_{0j} = 1 \quad j \in C \label{eq:om1_p0}
\end{align}

\noindent\emph{$u$-$p$ congruence}\\
These constraints impose the consistency between the ordering given by $u$ and $p$ variables.
\begin{align}
&u_i - u_j \ge 1 - (c+2)p_{ij} \quad i, j \in C, j\ne i\label{eq:om1_u_prec}\\
&u_i - u_j \le -1 + (c+2)(1 - p_{ij}) \quad i, j \in C, j\ne i\label{eq:om1_u_prec2}
\end{align}
\noindent\emph{Simple crossing sorties elimination}\\
Constraints \eqref{eq:on1_simpleCrossing} avoid that a sortie starts before a previous sortie is terminated. More specifically, if node $i$ is visited by the truck before node $l$ (i.e., $p_{il}=1$) and there are two sorties starting from $i$ and $l$, respectively, than the rendezvous node $k$ of the sortie from $i$, must be visited before node $l$.
\begin{align}
\begin{split}t_l^D \ge t_k^D - M
(3-\sum_{\substack{\tr{i}{j}{k}\in F\\j\neq l}}y_{ijk} - \sum_{\substack{\tr{l}{m}{n}\in F\\ m,n\not\in\{i,k\}}}y_{lmn}-p_{il})  \\
 i\in N_0, l \in C, k \in N_+,  i\neq l \neq k \end{split} \label{eq:on1_simpleCrossing}
\end{align}
\noindent \emph{Variable bounds}\\
\begin{align}
&t_i^T \ge 0 \quad i \in N\label{eq:om1_tt_bound}\\
&t_i^D \ge 0 \quad i \in N\label{eq:om1_td_bound}\\
&p_{ij}\in \left\{0,1\right\} \quad (i,j)\in A \label{eq:om1 p_bool}\\
&x_{ij}\in \left\{0,1\right\} \quad (i,j)\in A\label{eq:om1_x_bool}\\
&y_{ijk}\in \left\{0,1\right\} \quad \tr{i}{j}{k} \in F. \label{eq:om1_y_bool}\\
&1 \le u_i \le c+2 \quad i \in N_+\label{eq:om1_u_bound}
\end{align}
\subsection{Formulation ${\overline{\mbox{MC}}}$}\label{subsec:F1}
We propose a formulation called ${\overline{\mbox{MC}}}$  built upon MC with just a set of modifications inserted to provide a more realistic interpretation of the problem description given by Murray and Chu.

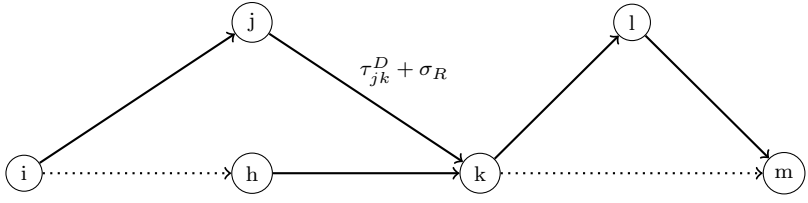
\begin{figure}
\centering
 \begin{tikzpicture}
 \node [draw, circle] at (-10,3) (i) {i} ;
 \node [draw, circle] at (-7,3) (h) {h} ;
 \node [draw, circle] at (-7,5) (j) {j} ;
 \node [draw, circle] at (-4,3) (k) {k};
 \node [draw, circle] at (-2,5) (l) {l};
 \node [draw, circle] at (0,3) (m) {m};
\draw[->, dotted,thick ] (i)--(h);
 \draw[->, thick] (i)--(j);
 \draw[->, thick] (h)--(k);
 \draw[->, thick] (j)--(k);
 \draw[->, thick] (k)--(l);
 \draw[->, thick] (l)--(m);
 \draw[->, dotted, thick] (k)--(m);
 \node at (-5,4.35) {$\tau_{jk}^D+\sigma_R$};
 \end{tikzpicture}
 \caption{Improving constraint.}
 \label{fig:const_wrong}
 \end{figure}
Let us consider  the \emph{truck timing constraints} \eqref{eq:om1_time5}, which applies to a truck running along arc $(h,k)\in A$.
Let us refer to Figure \ref{fig:const_wrong}, showing a solution where the drone performs a sortie which starts from $i$, a node that precedes $h$, terminates in $k$, and immediately starts a second sortie from $k$. If this happens, constraint \eqref{eq:om1_time5} imposes
$t_k^T \ge t_h^T + \tau_{hk}^T + \sigma_R+\sigma_L$, where $\sigma_R$ and  $\sigma_L$ refers, respectively, to the rendezvous and launch time arriving in $k$ and restarting from $k$.
Allocating the rendezvous time to $t_k^T$ is correct, but we believe that the launch time must be allocated to {the time of the node visited after $k$ in the truck route, otherwise the launching time needed for sortie $\langle k,l,m \rangle$ would be included in the flying time of sortie $\langle i,j,k \rangle$, and that could erroneously make exceed its endurance}. Moreover, we believe that the launch time must not be allocated to the truck when the drone starts  from the depot, since depot operations are considered to be done offline. Therefore,  we substituted \eqref{eq:om1_time5} with  \eqref{eq:f1_time5}-\eqref{eq:f1_time5b}.
\begin{align}
\begin{split}t_k^T \ge t_h^T + \tau_{hk}^T +  \sigma_R \sum_{\tr{i}{j}{k}\in F}  y_{ijk} +
\sigma_L\sum_{\substack{\tr{h}{r}{s}\in F\\ r\neq k}} y_{hrs} - M(1- x_{hk})  \\ \quad  (h,k)\in A, h\in N_+
\end{split} \label{eq:f1_time5}\\
\begin{split}
t_k^T \ge t_0^T + \tau_{0k}^T + \sigma_R(\sum_{\tr{i}{j}{k} \in F}y_{ijk} ) - M(1- x_{0k}) \quad (0,k)\in A, k\in N_+ \end{split} \label{eq:f1_time5b}
\end{align}

Moreover, we believe that $\sigma_L$ must  be included in the drone timing, when the drone is launched from a node different from the depot, being the time when the drone is on the truck and the operator prepares it.
We thus substitute the first \emph{drone timing constraint} \eqref{eq:om1_droneTime6}  with \eqref{eq:f1_droneTime6}--\eqref{eq:f1_droneTime6a}.
\begin{align}
&t_j^D \ge t_i^D + \tau_{ij}^D +\sigma_L - M (1-\sum_{\tr{i}{j}{k}\in F}y_{ijk} ) \quad (i,j)\in A, j \in C', i\in N_+\label{eq:f1_droneTime6}\\
&t_j^D \ge t_0^D + \tau_{0j}^D  - M (1-\sum_{\tr{0}{j}{k}\in F}y_{ijk} ) \quad j \in C'\label{eq:f1_droneTime6a}
\end{align}
\subsubsection{{Wait and no-wait models}}
As for the original MC formulation constraints \eqref{eq:f1_droneTime6}--\eqref{eq:f1_droneTime6a} allow $t^D_j$ to be greater than the arrival time at customer $j$,
without incurring in its computation while calculating energy consumption \eqref{eq:om1_energy}. {We refer to this model as the \emph{wait} model.}
From our experience the waiting is not always technically possible. In this case the  \emph{no-wait} model  is obtained by substituting the \emph{Drone battery endurance constraints} \eqref{eq:om1_energy} with
\eqref{eq:om1_energyRevised1}--\eqref{eq:om1_energyRevised2}.
\begin{align}
&t_k^D -t_i^D \le E + \sigma_L + M (1-y_{ijk} ) \quad \tr{i}{j}{k} \in F, i\neq 0 \label{eq:om1_energyRevised1}\\
&t_k^D -t_0^D \le E + M (1-y_{0jk} ) \quad \tr{0}{j}{k} \in F \label{eq:om1_energyRevised2}
\end{align}

Note  that our reformulation is similar, but slightly different from the one proposed by Ponza \cite{ponza2016optimization}.
Indeed, in \cite{ponza2016optimization} the launching time  is allocated to both truck and drone, so it enters in the computation of the energy consumption, while
we suppose that   launching time does not  consume energy, being the time needed by the truck to make the drone ready, while the drone does not fly.
Moreover, Ponza \cite{ponza2016optimization} includes the launching time in the computation also when the sortie starts from the depot, while we use the original MC model that assumes depot operations are done off-line.
\section{Improvements on DMN}\label{sec:ImprovementOnF1}
In this section we present a first enhanced formulation, built upon formulation ${\overline{\mbox{MC}}}$ of  Section \ref{subsec:F1}. The general idea is to substitute some explicit constraints with exponentially many constraints to be added in a cutting plane fashion.
\subsection{Crossing Sorties Elimination}
The basic MC model uses the explicit \emph{simple crossing sorties elimination} \eqref{eq:on1_simpleCrossing} which are based on timing and use of a `big M' constant. We propose to substitute them with the following structural  constraints that directly address and prevent the unfeasible  topologies.  Let $i \in N_0,l\in C$ be two {vertices visited} by the truck while running along a  path $P$ from $i$ to $l$, and assume that two sorties $\tr{i}{j}{k}$ and $\tr{l}{m}{n}$ with $k\not\in P$ exist. In this case, the second sortie starts before the first one is returned and the overall solution is therefore infeasible.
{Let us define $\mathcal{P}$ as the set of all the paths with these characteristics.}
The following  \emph{Crossing Sorties Elimination Constraints} (CSEC) \eqref{eq:om1_CSEC_1} can be used to eliminate this type of infeasibility.
\allowdisplaybreaks
\begin{align}
& \sum_{h=1}^{|P|-1} x_{v(h),v(h+1)} + \sum_{\substack{\tr{i}{j}{k}\in F\\k\not\in P}}y_{i,j,k} + \sum_{\tr{l}{m}{n}\in F}y_{l,m,n} \le |P| &P\in \mathcal{P}  \label{eq:om1_CSEC_1}
\end{align}
where $P=\{v(1),v(2),\dots, v(|P|)\}$, $v(1)=i, v(|P|)=l$.
These constraints impose that one of the arcs of the path or one of the sorties should be set to zero to make the solution feasible.

We can  strengthen \eqref{eq:om1_CSEC_1} using the idea of  ``tournament constraint". Since the path enters at most once in each vertex, we can include in the constraint also the arcs $(h,j)\in A$ such that $h,j \in P$ and $h$ precedes $j$ in $P$. We obtain the \emph{Tournament Crossing Constraints} (TCS), given by \eqref{eq:om1_CSEC_1}.
\allowdisplaybreaks
\begin{align}
& \sum_{h = 1}^{|P| - 1}\sum_{j  = h+1}^{|P|} x_{v(h)v(j)} + \sum_{\substack{\tr{i}{j}{k}\in F\\k\not\in P}}y_{i,j,k} + \sum_{\tr{l}{m}{n}\in F}y_{l,m,n}
 \le |P|&P\in \mathcal{P}\label{eq:om1_CSEC_2}
\end{align}

\subsection{Backward Sorties Elimination}

A next  improvement is the avoidance of backward sorties, that are sorties with a rendezvous that happens before the departure. Let $\mathcal{B}$ denote the set of all truck paths $P=\{v(1),v(2),\dots,v(q)\}$  with $v(1)=0$ and  $v(q)\in C$.  Given one of such paths $P$ suppose that exists a sorties $\tr{i}{j}{v(q)}$ with $i\not \in P$.
The {\em Backward Sortie Elimination Constraints} (BSEC) are given by \eqref{eq:om1_BSEC_1}, and impose that at least one of the arcs of the path or the sortie is eliminated.
\allowdisplaybreaks
\begin{align}
& \sum_{i = 1}^{|P| - 1} x_{v(i)v(i+1)} + \sum_{\tr{i}{j}{v(q)}\in F} y_{i,j,v(q)} \le |P| - 1 & P \in \mathcal{B}\label{eq:om1_BSEC_1}
\end{align}
Constraint \eqref{eq:om1_BSEC_1} can be strengthen by considering a tournament type constraint on path $P$, that imposes that at most one arc enters a node of the path,
and adding all sorties that terminate in $P$, but start at nodes outside $P$. We obtain the \emph{Tournament Backward Constraints} (TBS), given by \eqref{eq:om1_BSEC_2}.
\allowdisplaybreaks
\begin{align}
 &\mbox{(TBS)}
  \sum_{i = 1}^{|P| - 1}\sum_{j = i + 1}^{|P|} x_{v(i)v(j)} + \sum_{\substack{\tr{i}{j}{k} \in F\\i\not\in P, k\in P}} y_{i,j,k} \le  |P| - 1 &P\in \mathcal{B}
\label{eq:om1_BSEC_2}
\end{align}
Finally note that backwards sorties happen only when the time constraint are not respected, in our case in fractional solutions.
\subsection{Improving the Objective Function}
Preliminary computational results showed that formulation DMN, and thus F1, provides very bad lower bounds at the root node of the decision-tree (see Section \ref{sec:results} for more details).  Therefore, we reformulate the objective function $\min t^T_{c+1}$ by explicating its components:
the traveling time of the truck, the launching and rendezvous times, and the time truck waits for the drone.
The first two components can be modeled using the variables we dispose already, but the waiting times must be expressed by new variables. We  add variables $w_i, i \in N$, to model the truck waiting at node $i$,
thus giving the new objective function \eqref{eq:nm1_fobb}.
\allowdisplaybreaks
\begin{align}
\begin{split}
 \min  \sum_{(i,j)\in A}\tau_{ij}^Tx_{ij} + \sigma_R\sum_{\tr{0}{j}{k} \in F}y_{0jk} + (\sigma_L + \sigma_R) \sum_{\substack{\tr{i}{j}{k} \in F\\i\neq 0}}y_{ijk}  + \sum_{i \in N}w_i\end{split}\label{eq:nm1_fobb}
\end{align}
The  waiting variables must be included in the  truck timing constraints  as in \eqref{eq:nm1_truckTimeG}--\eqref{eq:nm1_timeLzero}. Note that we use a pair of constraint with opposite versus to impose equality  when  an arc $(h,k)$ is selected. Non-negativity of the  $w_i$ must also be imposed.
\begin{align}
\begin{split}t_k^T \ge t_h^T + \tau_{hk}^T +  \sigma_R \sum_{\tr{i}{j}{k}\in F}  y_{ijk} +
\sigma_L\sum_{\substack{\tr{h}{r}{s}\in F\\ r\neq k}} y_{hrs} \qquad\qquad\qquad\\- M(1- x_{hk})  + w_k\quad  (h,k)\in A, h\in N_+
\end{split} \label{eq:nm1_truckTimeG}\\
\begin{split}t_k^T \le t_h^T + \tau_{hk}^T +  \sigma_R \sum_{\tr{i}{j}{k}\in F}  y_{ijk} +
\sigma_L\sum_{\substack{\tr{h}{r}{s}\in F\\ r\neq k}} y_{hrs} \qquad\qquad\qquad\\+ M(1- x_{hk})  + w_k\quad  (h,k)\in A, h\in N_+
\end{split} \label{eq:nm1_truckTimeL}\\
\begin{split}
t_k^T \ge t_0^T + \tau_{0k}^T + \sigma_R(\sum_{\tr{i}{j}{k} \in F}y_{ijk} ) - M(1- x_{0k}) +w_k\quad (0,k)\in A, k \in N_+
\end{split} \label{eq:nm1_timeGzero}\\
\begin{split}
t_k^T \ge t_0^T + \tau_{0k}^T + \sigma_R(\sum_{\tr{i}{j}{k} \in F}y_{ijk} ) + M(1- x_{0k}) +w_k\quad (0,k)\in A, k \in N_+
\end{split} \label{eq:nm1_timeLzero}\\
\begin{split}w_i \ge 0 \quad i \in N \end{split} \label{eq:nm1_wait}
\end{align}
\subsection{{Subtour Elimination Constraints}}
{When avoiding crossing sorties by using the CSEC or TCS, Miller-Tucker-Zemlin constraint are not needed; indeed, total node ordering variables are not needed and timing constrains are sufficient to prevent subtours; however, the introduction of subtour elimination constraints \eqref{eq:om1_SEC_1} and the two nodes subtours elimination \eqref{eq:om1_SEC_2} is profitable. Preliminary computational results confirmed us the reduction of the number of visited nodes by the branch-and-bound and the computing time even if applied to the $\overline{\mbox{MD}}$ formulation. }

\allowdisplaybreaks
\begin{align}
\mbox{(SEC)}
&\sum_{i \in S}\sum_{j \in S} x_{ij} \le \left | S \right | -1&\hspace{-2cm}   & S \subseteq N, S\ne \emptyset, |S| > 2 \label{eq:om1_SEC_1}\\
&{x_{ij} + x_{ji} \leq 1} & & i, j \in N, i \neq j \label{eq:om1_SEC_2}
\end{align}\\

The overall formulation, called \textbf{DMN} in the following, uses all the improvements of this section with the ``tournament" version for crossing sorties \eqref{eq:om1_CSEC_2} and backward sorties \eqref{eq:om1_BSEC_2}.
\section{Two-indexed Formulations: DMN2}\label{sec:2form}
In the previous models we represented the sorties with three-indexed boolean variables, hereby we propose a new formulation which uses two-indexed arc variables. In particular, let us consider the binary variables $\overrightarrow{g}_{ij}$ and $\overleftarrow{g}_{jk}$ which take value 1 if the drone enters, respectively leaves,  customer node $j$. In this way the number of variables representing the sorties reduces from $n^3$ to $2n^2$. To further reduce the variables we preliminary fix to zero all the variables corresponding  to arcs with flying time exceeding the battery limit, i.e., we set $\overrightarrow{g}_{ij}=0$ for all $ (i,j)\in A: \tau^D_{ij}>E$ and
$\overleftarrow{g}_{jk}=0$  for all  $(j,k)\in A: \tau^D_{jk}+\sigma_R>E$. We also fix to zero variables that do not allow to complete a feasible drone fly: $\overrightarrow{g}_{ij}=0$, $(i,j)\in A, j\not\in C'$; $\overleftarrow{g}_{jk}=0$, $(j,k)\in A, j\not\in C'$; $\overrightarrow{g}_{ic+1}=0$ $i\in N$
 and $\overleftarrow{g}_{j0}=0$ $j\in N$.

Starting from our model DMN, we build the new formulation  DMN2 using the above two-indexed variables.
The objective function is
\begin{align}
 &\min  \sum_{(i,j)\in A}\tau_{ij}^Tx_{ij} + \sigma^L\sum_{\substack{(i,j)\in A\\i\neq 0}}\overrightarrow{g}_{ij} + \sigma^R\sum_{(j,k)\in A}\overleftarrow{g}_{jk} + \sum_{i \in N}w_i\label{eq:g2_fobb}
\end{align}
All the constraints using $y$ variables must be rewritten using the $\overrightarrow{g}$ and $\overleftarrow{g}$ variables, as follows:

\noindent
Customer covering  constraints \eqref{eq:om1_ass} become
\begin{align}
&\sum_{(i,j)\in A} x_{ij}  +  \sum_{(i,j)\in A} \overrightarrow{g}_{ij} = 1 \quad  j \in C\label{eq:f5_ass}\\
&\sum_{(i,j)\in A} x_{ij}  +  \sum_{(i,j)\in A} \overleftarrow{g}_{ij} = 1 \quad  i \in C\label{eq:f5_ass2}
\end{align}
Single sortie leave/return constrains \eqref{eq:om1_dep} and \eqref{eq:om1_arr} are no longer necessary, since are induced by the truck routing constraints \eqref{eq:om1_start}--\eqref{eq:om1_flow} within the $x$-$y$ coupling constraints, now  called
\noindent\emph{$x$-$g$ coupling constraints}:
\begin{align}
& \sum_{(i,j)\in A}\overrightarrow{g}_{ij} \le \sum_{(i,h)\in A}x_{ih} \quad i \in  N_0 \label{eq:f5_link}\\
& \sum_{(i,j)\in A}\overleftarrow{g}_{ij} \le \sum_{(h,j)\in A}x_{hj} \quad j \in  N_+ \label{eq:f5_link2}\\
& \sum_{(i,j)\in A} \overrightarrow{g}_{ij} = \sum_{(j,k)\in A} \overleftarrow{g}_{jk} \quad j \in C \label{eq:f5_flowg}
\end{align}
Drone-truck timing constraints \eqref{eq:om1_droneTruck1}--\eqref{eq:om1_droneTruck4} are substituted by
\begin{align}
&t_i^D \ge t_i^T - M (1-\sum_{(i,j)\in A}\overrightarrow{g}_{ij}) \quad i \in C\label{eq:f5_time1}\\
&t_i^D \le t_i^T + M (1-\sum_{i,j)\in A}\overrightarrow{g}_{ij}) \quad i \in C\label{eq:f5_time2}\\
&t_k^D \ge t_k^T - M (1-\sum_{(j,k)\in A}\overleftarrow{g}_{jk}) \quad k \in N_+\label{eq:f5_time3}\\
&t_k^D \le t_k^T + M (1-\sum_{(j,k)\in A}\overleftarrow{g}_{jk}) \quad k \in N_+\label{eq:f5_time4}
\end{align}
Truck timing constraints \eqref{eq:nm1_truckTimeG}--\eqref{eq:nm1_timeLzero} become
\begin{align}
\begin{split}t_k^T \ge t_h^T + \tau_{hk}^T + \sigma_R \sum_{(i,k)\in A}\overleftarrow{g}_{jk} +  \sigma_L\sum_{(h,l)\in A}\overrightarrow{g}_{hl} - M(1- x_{hk}) + w_k \\ \quad  (h,k)\in A, h\in N_+ \end{split}
\label{eq:g2_time5}\\
\begin{split}t_k^T \le t_h^T + \tau_{hk}^T + \sigma_R \sum_{(j,k)\in A}\overleftarrow{g}_{jk}  + \sigma_L\sum_{(h,l)\in A}\overrightarrow{g}_{hl} +  M(1- x_{hk}) + w_k \\ \quad  (h,k)\in A, h\in N_+ \end{split} \label{eq:g2_time5bis}\\
\begin{split}t_k^T \ge t_0^T + \tau_{0k}^T + \sigma_R \sum_{(j,k)\in A}\overleftarrow{g}_{jk}  - M(1- x_{0k}) + w_k  \quad k \in N_+ \end{split} \label{eq:g2_time5h0}\\
\begin{split}t_k^T \le t_0^T + \tau_{0k}^T + \sigma_R \sum_{(j,k)\in A}\overleftarrow{g}_{jk}  + M(1- x_{0k}) + w_k \quad  k \in N_+ \end{split} \label{eq:g2_time5bish0}\\
\end{align}
\noindent
Drone timing constraints \eqref{eq:f1_droneTime6}, \eqref{eq:f1_droneTime6a} and \eqref{eq:om1_droneTime7} are now\\
\begin{align}
&t_j^D \ge t_i^D + \tau_{ij}^D +\sigma_L - M(1-\overrightarrow{g}_{ij}) \quad j \in C',i \in C, i\neq j\label{eq:f5_time6}\\
&t_j^D \ge t_0^D + \tau_{0j}^D - M(1-\overrightarrow{g}_{0j}) \quad j \in C' \label{eq:f5_time6b}\\
&t_k^D \ge t_j^D + \tau_{jk}^D + \sigma_R - M (1-\sum_{(j,k)\in A}\overleftarrow{g}_{jk}) \quad j \in C',k\in N_+, j\neq k\label{eq:f5_time7}
\end{align}
Drone battery endurance constraints \eqref{eq:om1_energyRevised1}-\eqref{eq:om1_energyRevised2} become:
\begin{align}
\begin{split}t_k^D - t_i^D \le E + \sigma_L+ M (2 -\overrightarrow{g}_{ij} - \overleftarrow{g}_{jk}) \quad \\ (i,j)\in A, (j,k)\in A, j \in C',i \in N_+, i \neq k
\end{split}
\label{eq:f5_energy1}\\
\begin{split}t_k^D -t_0^D \le E + M (2 -\overrightarrow{g}_{0j} - \overleftarrow{g}_{jk}) \quad (j,k)\in A, j \in C'\end{split}
\label{eq:f5_energy2}
\end{align}
Crossing sorties elimination constraints \eqref{eq:om1_CSEC_2} must be rewritten as follows.
Let $i \in N_0,l\in C$ be two {vertices encountered} by the truck while running along a path $P$ from $i$ to $l$, and assume that exist two sorties defined by $\overrightarrow{g}_{ij}>0$ and $\overrightarrow{g}_{lm}>0$ and such that there is no node $k\in P\setminus\{i\}$ with $\overleftarrow{g}_{hk}>0$. In this case the second sortie starts before the first sortie is terminated an the following ``tournament" crossing sorties elimination holds:
\allowdisplaybreaks
\begin{align}
\mbox{(TCS2)}
& \sum_{h = 1}^{|P| - 1}\sum_{j  = h+1}^{|P|} x_{v(h)v(j)} + \sum_{\substack{(i,j)\in A, \\j\not\in P}}\overrightarrow{g}_{ij} +\sum_{\substack{(l,j)\in A, \\j\not\in P}}\overrightarrow{g}_{lj}
 \le |P|& P \in \mathcal{P}\label{eq:f5_CSEC_2}
\end{align}
where $P=\{v(1),v(2),\dots,v(q)\}$ with $v(1)=i, v(q)=l$, {and $\mathcal{P}$ now defines the set of all the paths with the described characteristics.} \\

\noindent Backward sorties elimination are modified as follows. Let $i \in N_0, j \in N_+$ be two nodes encountered by the truck while traveling on a path $P$, suppose that exist a sortie identified by $\overleftarrow{g}_{k,i}>0$ and such that  $k\not \in P$ and a sortie identified by $\overrightarrow{g}_{v(q),k}>0$ and such that  $k\not \in P$, that, together determine an infeasible solution. Let $\mathcal{B}$ now denote the set of all truck paths $\{v(1),\dots, v(m), \dots,v(q)\}$  with $v(1)=0$, $v(m) = i$, and  $v(q) = j$.
The tournament version of the backward sorties elimination constraints is:
\allowdisplaybreaks
\begin{align}
 &\mbox{(TBS2)}
  \sum_{h = 1}^{|P| - 1}\sum_{l = h + 1}^{|P|} x_{v(h)v(l)} + \overrightarrow{g}_{jk} +  \sum_{\substack{(l,k)\in A\\l\not\in P}} \overrightarrow{g}_{lk} + \sum_{\substack{(k,l)\in A\\l\not\in P}} \overleftarrow{g}_{kl} \le  |P| &P\in \mathcal{B}
\label{eq:f5_BSEC_2}
\end{align}

We finally  need to impose the following constraints to avoid infeasibilities:
\begin{align}
&\overrightarrow{g}_{ij} + \overleftarrow{g}_{ij} \le 1 \quad (i,j) \in A\\
&\overrightarrow{g}_{ij} + \overleftarrow{g}_{ji} \le 1 \quad (i,j) \in A
\end{align}

\section{Computational Experiments} \label{sec:results}
To test the above models we have implemented them to run on
 an Intel Core i3-2100 CPU, with 3.10 GHz and 8.00 GB of RAM, running Windows 7 operating system. CPLEX 12.71 was used as MILP solver,  and only a single thread was utilized during the testing.

 Formulation $\overline{\mbox{MC}}$ was solved directly by CPLEX, while DMN and DMN2 required a branch-and-cut implementation, since they include exponentially many constraints. To separate these  constraints we considered the residual graph $G'=(N,A')$ obtained from $G$ (see Section \ref{sec:prob}) by selecting the only arcs associated with a non zero variable ($x$, $y$, $\overrightarrow{g}$ and $\overleftarrow{g}$) in the continuous relaxation of the model. For the \emph{crossing sorties elimination} constraints and for the \emph{backward sorties elimination} constraints we simply explore the graph starting from depot 0, until a truck path violating one of the constraints is identified, if any. The overall procedure has a time complexity $O(|A'|)$. To separate the subtour elimination constraints we use the standard approach which requires to solve at most $O(n^2)$ max flow problems on the residual graph.

To facilitate the solution of the problem, we give all the algorithms an initial heuristic solution computed as follow.
We first find a heuristic (TSP) solution in which all the customers are served by the truck. We build  this solution with a greedy constructive algorithm followed by a local search improvement using as moves to define neighbouring solution the relocate of one customer and the swap of two customers. Next we examine all the possible sorties and we select the one, if any, that improves the solution. We repeat the search for improving sorties (and update of the solution) until no one is found.
We based our tests on the 72 benchmark instances provided by Murray and Chu \cite{murray2015flying}. In each of these, ten customers are randomly distributed across an eight-mile square region, while the depot is located in four different positions. The endurance of the drone was chosen to be either 20 or 40 minutes, while the speed of the drone was selected to be 15, 25, 35 miles/h based on Euclidean distances. The truck speed was assumed to be 25 miles/h and based on Manhattan distances. In the first 24 instances it has been set a ratio $|C'|/|C|=90\%$, while for the remaining ones they set the same ratio to 80\%. For each instance we tested the formulations on both the `no-wait' case (in which the drone is not allowed to wait at a customer), and the `wait' case, adopted  in \cite{murray2015flying}, where wait is allowed, but not considered in the computation of the battery endurance (see Section \ref{sec:intro}). Overall 144 instances are available in our test bed.

The interested reader can find all the results and solutions in our web site \url{www.or.unimore.it} following the \emph{online resources} link.
\subsection{Tests}

Table \ref{tab:lbroot} presents the performances of the lower bounds for different formulations.
The entries of the table are the average gap over the instances referred in each row, computed as  $100\cdot(opt-LB)/opt $, where `opt' is the value of the optimal solution obtained giving to formulation DMN2 enough time to close all instances.
The first column, labeled `E', reports the battery endurance value, while column labeled `speed' refers to the drone speed. Column  `$\overline{\mbox{MC}}$' reports the lower bound gap obtained by the Murray-Chu formulation, modified as in Section \ref{subsec:F1}, after one hour of computing time. All the other columns give the lower bound obtained at the root node.
Columns labeled `newF' (new objective function) refer to  formulation $\overline{\mbox{MC}}$ modified by adopting the objective function  introduced in Section \ref{sec:ImprovementOnF1}. Columns labeled `newF2' refer to formulation newF modified  by substituting the three index variables $y$ with the two-indexed variables $\overrightarrow{g}$ and $\overleftarrow{g}$ as in Section \ref{sec:2form}. Finally `DMN' and `DMN2' report on the  lower bound gaps of our complete formulations.
\begin{figure}[ht]
\caption{{Comparison of lower bound gaps at the root node for the different proposed methods.}}\label{fig:notsure}
\pgfplotsset{width=7cm,compat=1.8}
\begin{tikzpicture}
  \centering
  \begin{axis}[
        ybar, axis on top,
        title={No-wait instances},
        height=6cm, width=11.5cm,
        bar width=0.2cm,
        major grid style={draw=black},
        enlarge y limits={value=.1,upper},
        ymin=0, ymax=100,
        axis x line*=bottom,
        axis y line*=left,
        y axis line style={opacity=10},
        ytick={0, 50,  100}, tickwidth=5pt,
        enlarge x limits=true,
        legend style={
            at={(0.5,-0.5)},
            anchor=north,
            legend columns=-1,
            /tikz/every even column/.append style={column sep=0.6cm}
        },
        ylabel={Lower bound gap (root node)},
        symbolic x coords={E=20;sp=15, E=20;sp=25, E=20;sp=35, E=40;sp=15, E=40;sp=25, E=40;sp=35},
       xtick=data,
       x tick label style={rotate=45, anchor=east},
      nodes near coords={      \scriptsize
       \pgfmathprintnumber[precision=0]{\pgfplotspointmeta}
       }
    ]
    \addplot [draw=none, fill=blue!60] coordinates {
(E=20;sp=15,             83.31)
(E=20;sp=25,             81.99)
(E=20;sp=35,             79.72)
(E=40;sp=15,             82.86)
(E=40;sp=25,             80.68)
(E=40;sp=35,             79.66)
 };
    \addplot [draw=none, fill=red!60] coordinates {
(E=20;sp=15,             14.38)
(E=20;sp=25,             26.60)
(E=20;sp=35,             19.47)
(E=40;sp=15,             29.95)
(E=40;sp=25,             23.58)
(E=40;sp=35,             20.93)
 };
     \addplot [draw=none, fill=green!70] coordinates {
(E=20;sp=15,             13.71)
(E=20;sp=25,             25.30)
(E=20;sp=35,             19.14)
(E=40;sp=15,             29.27)
(E=40;sp=25,             21.94)
(E=40;sp=35,             19.85)
 };
     \addplot [draw=none, fill=blue!80] coordinates {
(E=20;sp=15,             9.35)
(E=20;sp=25,             23.10)
(E=20;sp=35,             16.12)
(E=40;sp=15,             26.53)
(E=40;sp=25,             20.05)
(E=40;sp=35,             18.22)
 };
    \addplot [draw=none, fill=red!80] coordinates {
(E=20;sp=15,             8.46)
(E=20;sp=25,             18.57)
(E=20;sp=35,             11.41)
(E=40;sp=15,             22.90)
(E=40;sp=25,             16.06)
(E=40;sp=35,             12.60)
 };

    \legend{\small $\overline{\mbox{MC}}$, newF, newF2, DMN, DMN2}
  \end{axis}
\end{tikzpicture}
\end{figure}
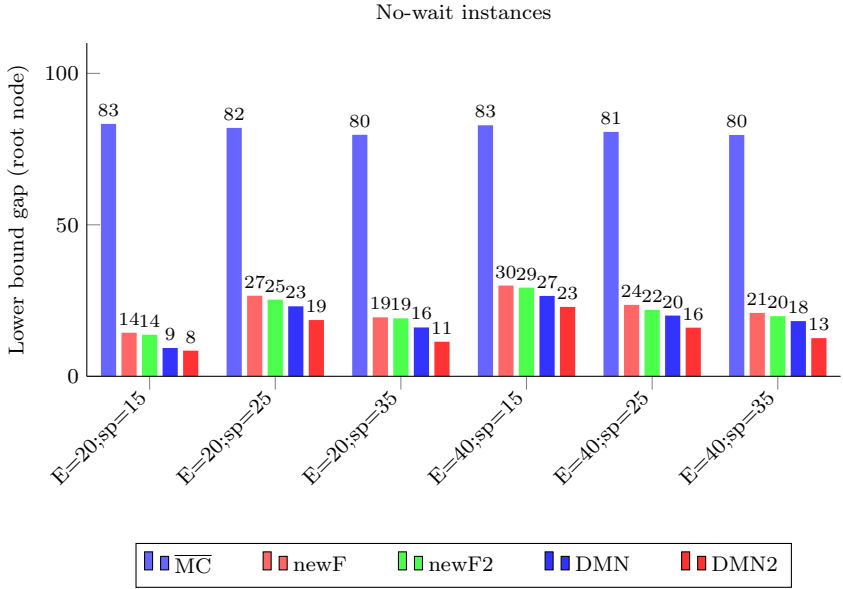
\setlength\tabcolsep{3pt}
\begin{table}[ht]
\scriptsize
  \centering
  \caption{Lower bounds at the root node, by speed class} \label{tab:lbroot}
    \begin{tabular}{lcrrrrrrrrrrrrrrrrr}
    \toprule
    &&&\multicolumn{5}{c}{No-Wait}&&\multicolumn{5}{c}{Wait}\\
$E$ & speed&& $\overline{\mbox{MC}}$ &newF & newF2 &DMN&DMN2&&$\overline{\mbox{MC}}$ &newF & newF2 &DMN&DMN2    \\
\cline{1-2}\cline{4-8}\cline{10-14}
     &15    &&   83.31        &14.38     &13.71     &9.35      &8.46      &&83.42     &13.83     &12.73     &8.67      &6.83   \\
\bf20&25    &&   81.99        &26.60     &25.30     &23.10     &18.57     &&81.39     &23.99     &21.79     &20.25     &15.36  \\
     &35    &&   79.72        &19.47     &19.14     &16.12     &11.41     &&79.64     &19.21     &18.56     &15.99     &11.12  \\
\cline{1-2}\cline{4-8}\cline{10-14}
     &15    &&   82.86        &29.95     &29.27     &26.53     &22.90     &&82.78     &29.77     &29.30     &27.16     &23.06  \\
\bf40&25    &&   80.68        &23.58     &21.94     &20.05     &16.06     &&80.68     &23.58     &22.55     &20.28     &15.80  \\
     &35    &&   79.66        &20.93     &19.85     &18.22     &12.60     &&79.65     &20.79     &19.61     &17.82     &12.49  \\
\cline{1-2}\cline{4-8}\cline{10-14}
\multicolumn{2}{c}{\bf All}
            &&\bf81.37        &\bf22.49  &\bf21.53  &\bf18.89  &\bf15.00  &&\bf81.26  &\bf21.86  &\bf20.76  &\bf18.36  &\bf14.11\\
\end{tabular}
\end{table}
Formulation $\overline{\mbox{MC}}$ clearly gives bad results, in terms of lower bound gap, providing a gap close to 80\% for all class of instances.
The adoption of the new objective function greatly improves the behaviour, drastically reducing the gap to about 22\%. The use of two-indexed variables and complete formulations further improve the gap. Concerning the effect of the waiting time at the customers we observe that instances where is allowed to wait have slightly smaller gaps.

Table \ref{tab:lbrootDepot} gives the same data as in Table \ref{tab:lbroot}, but grouped by depot position. Column labeled ‘dep’ refers to the four depot positions, with respect to the square in which are generated the customers. In the first position ‘a’ the depot is almost in the barycentre of the customers, in ‘b’, ‘c’ and ‘d’ it is close to the right border of the square, with a vertical position which is, respectively, in the barycentre, on the bottom border and below the bottom border at a distance which is equal to the distance of the barycenter from this border. One can see that the worst bound gaps are obtained when the depot is in position `b' and improve when it is far from the barycenter (positions `c' and `d'). The best gaps are obtained when the depot is outside the square, that is far from all customers. In this case there is a sort of `offset' distance that must be covered by any route to reach the customers, and this improves the bound.

\setlength\tabcolsep{3pt}
\begin{table}[ht]
\scriptsize
  \centering
  \caption{Lower bounds at the root node, by depot location} \label{tab:lbrootDepot}
    \begin{tabular}{lcrrrrrrrrrrrrrrrrr}
    \toprule
    &&&\multicolumn{5}{c}{No-Wait}&&\multicolumn{5}{c}{Wait}\\
$E$ & dep&& $\overline{\mbox{MC}}$ &newF & newF2 &DMN&DMN2&&$\overline{\mbox{MC}}$ &newF & newF2 &DMN&DMN2    \\
\cline{1-2}\cline{4-8}\cline{10-14}
     & a &&89.87      &20.97     &20.24     &15.88     &13.79     &&89.67     &19.10     &17.17     &13.80     &11.21\\
\bf20& b &&88.08      &25.19     &24.25     &21.57     &16.28     &&88.14     &24.16     &22.81     &19.86     &14.72\\
     & c &&82.06      &19.03     &17.86     &14.53     &11.10     &&81.94     &18.46     &17.67     &14.91     &9.95 \\
     & d &&66.69      &15.41     &15.18     &12.78     &10.08     &&66.18     &14.31     &13.11     &11.31     &8.55 \\
\cline{1-2}\cline{4-8}\cline{10-14}
     & a &&90.08      &26.65     &25.03     &23.00     &18.71     &&90.07     &26.46     &24.60     &23.07     &18.84\\
\bf40& b &&87.41      &29.92     &28.72     &26.86     &21.08     &&87.41     &29.92     &28.72     &27.55     &20.75\\
     & c &&81.30      &23.49     &21.88     &20.16     &16.01     &&81.29     &23.46     &23.42     &19.87     &15.79\\
     & d &&65.46      &19.21     &19.11     &16.38     &12.95     &&65.37     &19.01     &18.53     &16.51     &13.08\\
     \bottomrule
\end{tabular}
\end{table}
 Figure \ref{fig:notsure} provides a graphical representation of the lower bound gaps for the 'no-wait' instances.
{
One can see that the lower bound gap at the root node is improved largely by the proposed methods with respect to the original model with no difference based on endurance and speed class. The reader can notice that each of the improvement features applied provide a relevant contribution to the results.
}

In Tables \ref{tab:exactWaitSpeed}-\ref{tab:exactWaitDepot} we report on the performances of the formulations in finding the proven optimal solution for the instances where waiting at the customer's location is allowed. The values in the table are averaged over the 12 (resp. 9) instances of each class of speed (resp. depot location). The overall results are computed over the 72 instances.
In the columns labeled `gap\%', `opt', `nodes' and `time', we report, respectively: the percentage gap computed as $100\cdot (UB-LB)/UB$ (average on the only instances not solved to the optimum), the number of proven optimal solutions, the average computing time in CPU seconds, and the average number of branch-decision-trees explored.
As already anticipated, formulation $\overline{\mbox{MC}}$, in one hour of CPU time,  was not able to solve any instance given the poor performances of the lower bound, so we do not report the CPU time for this formulation.

Looking at Table \ref{tab:exactWaitSpeed} we observe that the instances with endurance 20 are much easier than those with endurance 40: in terms of obtained optimal solutions and computing times. This is due to the fact that the number of possible drone sorties are much less in the first case, thus reducing the possible choices of the algorithm. The two-indexed formulation DMN2 has better performances than DMN, since it solves 13 more instances, the unsolved instances have smaller gap and also the computing time is smaller.

For the 'no-wait' instances the formulations exhibit a very similar behaviour that is not reported here {(the interested reader can find the related results at  \url{www.or.unimore.it})}.

Looking at Table \ref{tab:exactWaitDepot} we can observe that instances of class `a', where the depot is located almost in the middle of the customers, are more difficult than the other cases. Class `d', which gave the best results for the root node lower bounds, has, instead, the same difficulty of the classes `b' and `c'

\setlength\tabcolsep{2pt}
\begin{table}[ht]
  \centering
  \scriptsize
  \caption{Exact solutions for `wait' instances, by speed class} \label{tab:exactWaitSpeed}
    \begin{tabular}{lcrrrrrrrrrrrrrrrrr}
    \toprule
          & && \multicolumn{3}{c}{$\overline{\mbox{MC}}$} &&\multicolumn{4}{c}{DMN}&&\multicolumn{4}{c}{DMN2}\\
\cline{4-6}\cline{8-11}\cline{13-16}
$E$ & sp&& gap\%&opt&nodes&& gap\%&opt&time&nodes&&gap\%&opt&time&nodes&\\
\cline{1-2}\cline{4-6}\cline{8-11}\cline{13-16}
  &15        &&83.48     &0         &180356.6  &&0.00      &12        &10.2     &379.9    && 0.00     & 12        &6.1      &340.4     \\
20&25        &&82.02     &0         &159165.3  &&8.25      &11        &797.1    &17056.3  && 1.75     & 11        &503.1    &19934.1   \\
  &35        &&80.81     &0         &142095.8  &&8.72      &9         &1342.2   &24671.7  && 4.18     & 10        &920.4    &35803.4   \\
\cline{1-2}\cline{4-6}\cline{8-11}\cline{13-16}
  &15        &&82.96     &0         &153948.3  &&10.59     &2         &3313.3   &54078.5  && 7.67     & 9         &2247.2   &91281.3   \\
40&25        &&81.47     &0         &145158.7  &&10.14     &4         &2533.3   &46797.5  && 6.70     & 7         &2033.7   &83462.6   \\
  &35        &&80.63     &0         &140049.8  &&5.28      &8         &1433.6   &25229.2  && 8.93     & 10        &1020.6   &45452.8   \\
\cline{1-2}\cline{4-6}\cline{8-11}\cline{13-16}
\multicolumn{2}{c}{\bf All}  &&\bf81.90 &\bf0     &\bf153462.4        &&\bf9.33 &\bf46        &\bf1571.68         &\bf28035.5         &&\bf6.50 &\bf59        &\bf1121.89         &\bf46045.8\\
     \bottomrule
\end{tabular}
\end{table}
\setlength\tabcolsep{2pt}
\begin{table}[ht]
  \centering
  \scriptsize
  \caption{Exact solutions for `wait' instances, by depot location} \label{tab:exactWaitDepot}
    \begin{tabular}{lcrrrrrrrrrrrrrrrrr}
    \toprule
          & && \multicolumn{3}{c}{$\overline{\mbox{MC}}$} &&\multicolumn{4}{c}{DMN}&&\multicolumn{4}{c}{DMN2}\\
\cline{4-6}\cline{8-11}\cline{13-16}
$E$ & dep&& gap\%&opt&nodes&& gap\%&opt&time&nodes&&gap\%&opt&time&nodes&\\
\cline{1-2}\cline{4-6}\cline{8-11}\cline{13-16}
  &a         &&90.14     &0         &161330.4  &&9.03      &6         &1496.6    &26211.3   &&3.37      &6         &1293.5    &41864.8\\
20&b         &&88.38     &0         &162479.0  &&7.34      &8         &1034.0    &22432.4   &&0.00      &9         &559.4     &29790.4\\
  &c         &&82.72     &0         &157777.2  &&0.00      &9         &155.0     &3771.9    &&0.00      &9         &29.2      &1703.1 \\
  &d         &&67.18     &0         &160570.2  &&0.00      &9         &180.7     &3728.2    &&0.00      &9         &24.2      &1412.2 \\
\cline{1-2}\cline{4-6}\cline{8-11}\cline{13-16}
  &a         &&90.64     &0         &145540.9  &&11.71     &2         &2937.7    &49395.2   &&8.23      &3         &2662.4    &86080.9\\
40&b         &&87.79     &0         &145526.4  &&10.68     &3         &2555.8    &46533.9   &&7.67      &6         &2081.0    &86932.6\\
  &c         &&82.02     &0         &151308.9  &&7.51      &5         &2033.6    &35907.6   &&0.00      &9         &1064.3    &53117.3\\
  &d         &&66.30     &0         &143166.2  &&6.41      &4         &2180.0    &36303.6   &&1.96      &8         &1261.1    &67464.7\\
\cline{1-2}\cline{4-6}\cline{8-11}\cline{13-16}
\multicolumn{2}{c}{\bf All}  &&\bf81.90	&\bf0	&\bf153462.4	&&\bf9.33	&\bf46	&\bf1571.7	&\bf28035.5	&&\bf6.50	&\bf59	&\bf1121.9	&\bf46045.8\\
     \bottomrule
\end{tabular}
\end{table}

{Tables \ref{tab:waitnowaitspeeddepot} depicts the difference between the {'wait'} and {'no-wait'} solutions with respect to the speed class and the depot location. The entries of Table \ref{tab:waitnowaitspeeddepot}-(a) and Table  \ref{tab:waitnowaitspeeddepot}-(b), are the endurance, the speed (depot location, respectively), the gap$\%$ computed as $100\cdot (opt^{n}-opt^{w})/opt^{n}$ averaged only when the two solution differs, being {$opt^{n}$ and $opt^{w}$ the optimal solution of the no-wait and wait case, respectively.} The last column reports  the number of occurrences of a difference between the two optimal solutions.  One can note that a more restrictive endurance make the allowance of waiting at customer without flying more important. This appears more relevant than the speed class. On the other hand, depot position 'a' (centrally positioned) and 'd' (positioned out of the square) have a strong impact on the solution when waiting is allowed. Moreover, based on extensive computational tests, we can state that the 'wait' model is only slightly easier to be solved than the 'no-wait' one.
}

\begin{figure}[h]
    \caption{{Comparison between the 'no-wait' (a) and 'wait' (b) solution of an instance 
    with $E=20$, depot location 'd', and 'speed' = 25 miles/h. The solid, dashed, lines represent the truck route, sorties, respectively. The square indicates the depot while the circles indicate the customers.}}
     \centering
     \subfigure[no-wait]{
         \centering
         \includegraphics[width=0.4\textwidth]{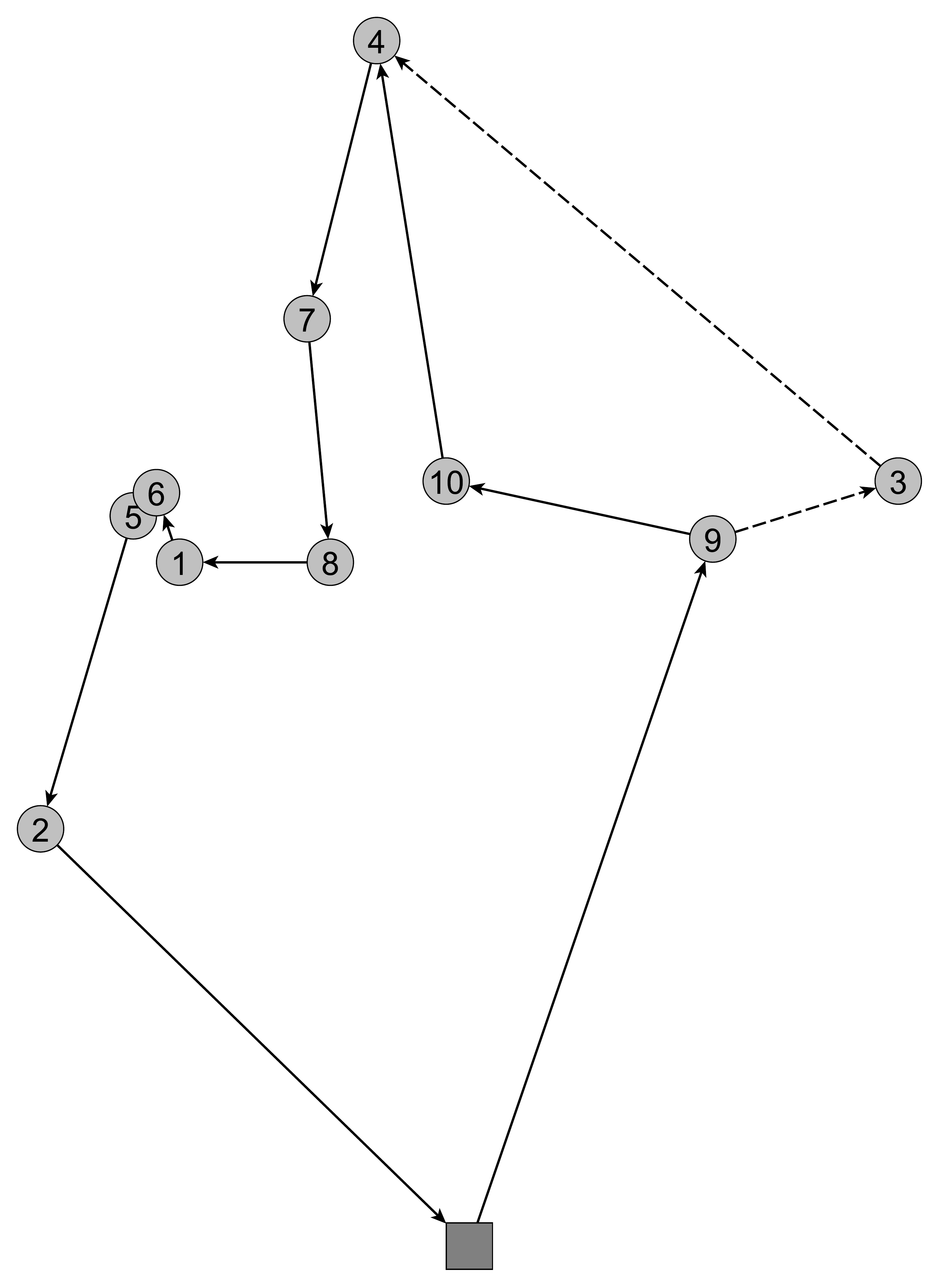}
         \label{fig:nowait}
     }
     \hfill
     \subfigure[wait]{
         \centering
         \includegraphics[width=0.4\textwidth]{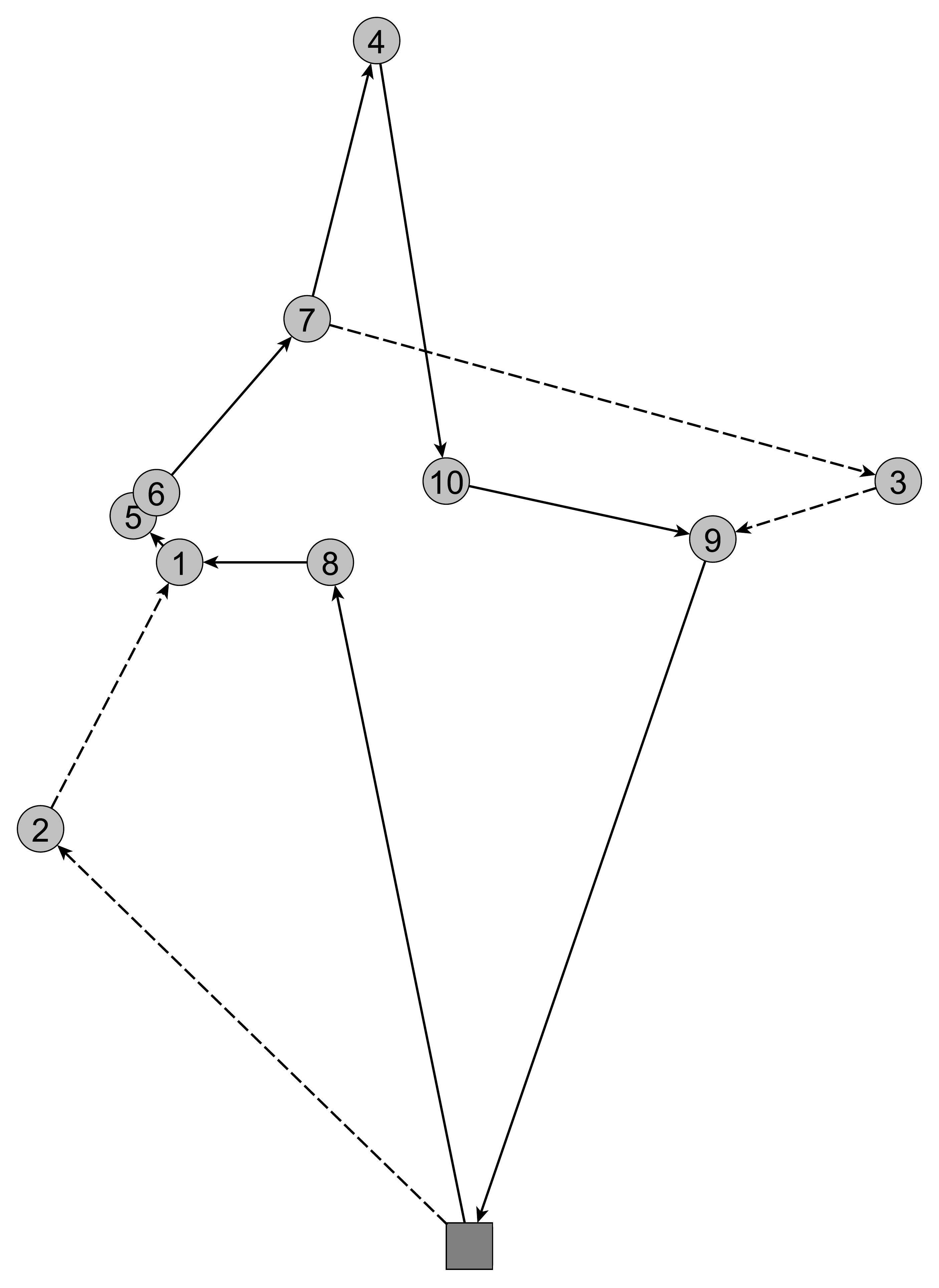}
         \label{fig:wait}
     }
        \label{fig:waitnowait}
\end{figure}

\begin{table}[h]
\scriptsize
\centering
	\caption{Comparison between {'wait'} and {'no-wait'} solutions}
	\subtable[By speed class]{
		\centering
		    \begin{tabular}{rrrr}
    \toprule
    \multicolumn{1}{l}{$E$} & \multicolumn{1}{l}{speed} & \multicolumn{1}{l}{gap$\%$} & \multicolumn{1}{l}{occur.}\\
   \cmidrule(lr){1-2} \cmidrule(lr){3-4}
    20 & 15 & 1.65 & 5 \\
      & 25 & 3.86 & 10 \\
      & 35 & 2.24 & 2 \\
   \cmidrule(lr){1-2} \cmidrule(lr){3-4}
    40 & 15 & 1.42 & 4 \\
      & 25 & 0.00 & 0 \\
      & 35 & 0.00 & 0 \\
   \cmidrule(lr){1-2} \cmidrule(lr){3-4}
     All &  & 2.72 & 21 \\
      \bottomrule
    \end{tabular}%
	}
	\hspace{1cm}
	\subtable[By depot location]{
		\centering
	    \begin{tabular}{rrrr}
    \toprule
    \multicolumn{1}{l}{$E$} & \multicolumn{1}{l}{dep} & \multicolumn{1}{l}{gap$\%$} & \multicolumn{1}{l}{occur.} \\
      \cmidrule(lr){1-2} \cmidrule(lr){3-4}
    20 & \multicolumn{1}{l}{a} & 4.08 & 5 \\
      & \multicolumn{1}{l}{b} & 2.64 & 5 \\
      & \multicolumn{1}{l}{c} & 1.22 & 4 \\
      & \multicolumn{1}{l}{d} & 3.89 & 3 \\
        \cmidrule(lr){1-2} \cmidrule(lr){3-4}
    40 & \multicolumn{1}{l}{a} & 2.66 & 1 \\
      & \multicolumn{1}{l}{b} & 0.00 & 0 \\
      & \multicolumn{1}{l}{c} & 0.37 & 1 \\
      & \multicolumn{1}{l}{d} & 1.32 & 2 \\
        \cmidrule(lr){1-2} \cmidrule(lr){3-4}
     All &   & 2.72 & 21 \\
     \bottomrule
    \end{tabular}%
	}
	\label{tab:waitnowaitspeeddepot}
\end{table}
{Figure \ref{fig:waitnowait} shows the difference between a 'no-wait' and a 'wait' solution. For this particular instance, the gap between the value of the two optimal solutions is around 9$\%$. Note that one more sortie can be performed in the 'wait' solution. In particular, allowing to wait at customer 2, the drone can save battery while the truck travels path (0,8,1), whose length would exceed battery endurance if the drone had to wait while flying. Similarly, the 'wait' solution can let the drone wait at customer 3, while the truck runs on the path (7,4,10,9), whose time would exceed the battery endurance. In this case the 'no-wait' solution performs a shorter sortie, instead; however, this difference does not improve the solution, but it shows, nevertheless, that the instances of the FSTSP can have more than one optimal solution, with the same truck route, and different sorties. This last consideration is one of the reasons why these problems are hard to solve also for small instances and it further justifies the increasing interest of researchers.}
\section{Conclusions}\label{sec:conclusions}
{This paper considered one of the most promising fields in parcel deliveries nowadays, the combined use of traditional vehicles and drones.
We focused on one particular problem called the Flying Sidekick Traveling Salesman Problem, where a drone and a truck are coupled and synchronized, and must serve a set of customers.
We provided improvements on the literature by proposing a three and a two-indexed formulations for which we proposed a novel objective function that could increase remarkably the lower bounds, and a set of good inequalities to be separated in a branch-and-cut fashion that provided an important contribution to obtain better solutions in a faster way.
 Our method outperforms the previously proposed one, being the two-indexed formulation the preferable one: 59 out of the 72 benchmark instances could be solved to optimality with an average percentage gap between the best lower and upper bound of 1.70\% and within an average computing time of less than 20 minutes.
 We evaluated two versions of the problem: one in which the drone is allowed to wait at the customers to save battery and one where this is prohibited. Both are hard to solve even when considering small sized instances.
 The version in which waiting is allowed is only slightly easier to be solved with respect to the other one: larger sets of feasible solutions must be treated by the algorithm, but feasible solutions are easier to find and possibly with a smaller cost. Among the several features considered, the drone endurance is the one that has the stronger influence on the convergence of the algorithm, a smaller endurance allow the algorithms to have better performances, while reducing the number of feasible sorties.
To conclude, our ideas could be easily applied to similar problems: this could suggest a direction of future research where analogous problems, with different constraints, are considered. Another possible direction of research could be the development of metaheuristic algorithms for these sets of problems.}


\bibliographystyle{spmpsci}      
\bibliography{references}   


\end{document}